\documentclass[doublespace]{article}

\usepackage{amssymb}
\usepackage{amsmath}   % For Latex2e
\usepackage{amsfonts}
\usepackage{amsthm}

\newtheorem{teorema}{Theorem}[section]
\newtheorem{lema}[teorema]{Lemma}
\newtheorem{proposicion}[teorema]{Proposition}

\newtheorem{conjetura}[teorema]{Conjecture}

\theoremstyle{definition}

\newtheorem{ejemplo}[teorema]{Example}

\newtheorem{contraejemplo}[teorema]{Counterexample}

\theoremstyle{remark}
\newtheorem{nota}[teorema]{Remark}

\newcommand{\code}{\mathcal{C}}

\newcommand{\fq}{\mathbb{F}_q}

\newcommand{\lgot}{\mathcal{L}}

\newcommand{\ze}{\mathbb{Z}}

\newcommand{\rr}{\mathbb{R}}
\newcommand{\mr}{M_{\mathbb{R}}}
\newcommand{\nr}{N_{\mathbb{R}}}
\newcommand{\zz}{\mathbb{Z}}

\newcommand{\codh}{\mathrm{H}^0 (X_P,\mathcal{O}(D_P))}

\newcommand{\talque}{~ | ~}
\newcommand{\ev}{\mathrm{ev}}
\renewcommand{\div}{\mathrm{div} }
\newcommand{\volr}{\mathrm{vol}_r}
\newcommand{\vold}{\mathrm{vol}_2}
\newcommand{\volt}{\mathrm{vol}_3}

\newenvironment{pf}{\noindent{\it Proof}.}{\qed}

\title{On the Parameters of $r$-dimensional Toric Codes}

\author{Diego Ruano\footnote{Partially supported by MEC MTM2004-00958 (Spain). Address: Department of Algebra, Geometry and Topology, Faculty of Sciences, University of Valladolid, E-47005 Valladolid, Spain. E-mail: ruano@agt.uva.es}}

\begin{document}

\maketitle

\begin{abstract}
From a rational convex polytope of dimension $r\ge 2$ J.P. Hansen
constructed an error correcting code of length $n=(q-1)^r$ over
the finite field $\fq$. A rational convex polytope is the same
datum as a normal toric variety and a Cartier divisor. The code is
obtained evaluating rational functions of the toric variety
defined by the polytope at the algebraic torus, and it is an
evaluation code in the sense of Goppa. We compute the dimension of
the code using cohomology. The minimum distance is estimated using
intersection theory and mixed volumes, extending the methods of
J.P. Hansen for plane polytopes. Finally we give a counterexample
to Joyner's conjectures \cite{jo}.
\end{abstract}

\section{Introduction}

An important family of error correcting codes is the
Algebraic-Geometry Codes, introduced by Goppa in 1981. These codes
became important in 1982, when Tsfasman, Vl\u adu\c t and Zink
constructed a sequence of error correcting codes that exceeds the
Gilbert-Varshamov bound. This was the first improvement of that
bound in thirty years.

The Algebraic-Geometry codes are defined by evaluating rational
functions on a smooth projective curve over a finite field $\fq$.
The functions of $\lgot (D)$ are evaluated in certain rational
points of the curve, where $D$ is a divisor whose support does not
contain any of the rational points where we evaluate. Their
parameters are estimated easily using the Riemann-Roch theorem
because the points can be seen as divisors.

This construction can be extended to define codes using normal
varieties of any dimension \cite{tsvl} giving rise to the called
evaluation codes. One can evaluate rational functions but the
estimation of the parameters is not easy in general, in particular
the estimation of the minimum distance is difficult.

The toric geometry studies varieties that contain an algebraic
torus as a dense subset and furthermore the torus acts on the
variety. The importance of these varieties, called toric
varieties, is based on their correspondence with combinatorial
objects, this makes the techniques to study the varieties (such as
cohomology theory, intersection theory, resolution of
singularities, etc) more precise and the calculus easier.

J.P. Hansen in 1998 (see \cite{ha,ha3}) considered evaluation
codes defined over some toric surfaces, in order to use the proper
combinatorial techniques of toric surfaces to estimate the
parameters of these codes. D. Joyner in 2004 (see \cite{jo}) also
considered toric codes over toric surfaces and he gave examples
with good parameters using a library in Magma to compute them. He
also proposed several questions and conjectures. Recently, other
works on toric codes have been published \cite{li,li2}.

This work treats evaluation codes over a toric varieties of
arbitrary dimension ($r \ge 2$) and length $(q-1)^r$ over the
finite field of $q$ elements. A rational convex polytope, is the
same datum as a normal toric variety and a divisor. For each
rational convex polytope we define an evaluation code over its
associated toric variety. The dimension of the code is computed
using cohomology theory, by the computation of the kernel of the
evaluation map. The minimum distance is estimated using
intersection theory and mixed volumes. Finally we give a
counterexample to the two conjectures of Joyner \cite{jo}.

We mainly use the notation of \cite{fu} for toric geometry
concepts and for all the toric geometry concepts and results we
refer to \cite{fu} and \cite{od}.

\section{Toric Geometry}

Let $N$ be a lattice ($N \simeq \zz^r$ for some $r\ge 1$). Let $M
= \mathrm{Hom} (N, \ze)$ be the dual lattice of $N$. One has the
dual pairing $\langle ~ , ~ \rangle: M \times N \to \zz$, $(u,v)
\mapsto u(v)$, that is $\zz$-bilinear. Let $\nr = N \otimes \rr$
and let $\mr = M \otimes \rr$. $\mr$ is the dual vector space of
$\nr$. One has the dual pairing $\langle ~ , ~ \rangle: \mr \times
\nr \to \mathbb{R}$, $(u,v) \mapsto u(v)$, that is
$\mathbb{R}$-bilinear.

Let $\fq$ be the finite field of $q$ elements and $T=(\fq^\ast)^r$
the $r$-dimensional \textbf{algebraic torus}. Let $\sigma$ be a
strongly convex rational cone in $\nr$ ($\sigma \cap (-\sigma) =
\{ 0 \}$ and $\sigma$ is generated by vectors in the lattice), for
the sake of simplicity we will just use the word \textbf{cone} in
this work. And let $\sigma^\vee$ be its dual cone $ \sigma^\vee =
\{ u \in \mr \talque \langle u,v\rangle \ge 0 ~ ~ \forall ~ v \in
\sigma\} $. A \textbf{face} $\tau$ of $\sigma$ is the intersection
with any supporting hyperplane.

Let $\sigma$ be a cone, then $S_\sigma = \sigma^\vee \cap M$ is a
finitely generated semigroup by Gordan's lemma. We consider its
associated $\fq$-algebra, $\fq [S_\sigma] = \bigoplus_{u \in
S_\sigma} \fq \chi^u$ ($\chi^u \chi^{u'} = \chi^{u+u'}$, the unit
is $\chi^0$). Therefore one can consider $U_\sigma = \mathrm{Spec}
(\fq [S_\sigma])$ which is the \textbf{toric affine variety
associated to $\sigma$}.

One can consider $\chi^u$ as Laurent monomial, $\chi^u (t) =
t_1^{u_1} \cdots t_r^{u_r} \in \fq[t_1, \ldots, t_r]_{t_1\cdots
t_r}$, this also gives a function $T \to \fq^\ast$. In the theory
of algebraic groups this is called a \textbf{character}.

A \textbf{fan} $\triangle$ in $N$ is a finite set of cones in
$\nr$ such that: Each face of a cone in $\triangle$ is also a cone
in $\triangle$ and the intersection of two cones in $\triangle$ is
a face of each. For a fan $\triangle$ the \textbf{toric variety}
$X_\triangle$ is constructed taking the disjoint union of the
affine toric varieties $U_\sigma$ for $\sigma \in \triangle$, and
gluing the affine varieties of common faces.

A toric variety is a disjoint union of orbits by the action of the
torus $T$. There is a one to one correspondence between
$\triangle$ and the orbits. For a cone $\sigma$ we denote by
$V(\sigma)$ the closure of the the orbit of $\sigma$, and one has
that $\dim \sigma + \dim V(\sigma) = r$.

A toric variety defined from a fan $\triangle$ is non-singular if
and only if for each $\sigma \in \triangle$, $\sigma$ is generated
by a subset of a basis of $N$. We say that a fan $\triangle'$ is a
\textbf{refinement} of $\triangle$  if each cone of $\triangle$ is
union of cones in $\triangle'$. One has a morphism $X(\triangle')
\to X (\triangle)$ that is birational and proper. By refining a
fan we can resolve the singularities considering a non-singular
refined fan, we assume in this work that a fan is always refined
and therefore its associated toric variety is non-singular.

A convex rational polytope in $\mr$ is the convex hull of a finite
set of points in $M$, for the sake of simplicity we just say
\textbf{polytope}. One can represent a polytope as the
intersection of halfspaces. For each facet $F$ (face of
codimension 1) there exists $v_F \in N$ inward and primitive and
an integer $a_F$ such that
$$
P = \bigcap_{F~\mathrm{is~a~facet}} \{ u \in \mr \talque  \langle
u, v_F \rangle \ge -a_F \}
$$

Given a face $p$ of $P$, let $\sigma_p$ be the cone generated by
$v_F$ for all the facets $F$ containing $p$. Then
$$
\triangle_P = \{ \sigma_p \talque p~\mathrm{is~a~face~of~} P  \}
$$ is a fan which is called \textbf{fan associated to $P$} and its associated toric
variety is denoted by $X_P$. We assume that the associated fan is
non-singular, in other cases we refine the fan and therefore we
consider the halfspaces associated to the new borders (see
\cite[section 5.4]{ge}).

From a polytope one can define the following $T$-invariant Weil
divisor (which is also a Cartier divisor because the variety is
non-singular),
$$
D_P = \sum_{F~\mathrm{is~a~facet}} a_F V(\rho_F)
$$ and given $u \in P$
$$
\mathrm{div}(\chi^u) = \sum_{F~\mathrm{is~a~facet}} \langle u, v_F
\rangle V(\rho_F)
$$

We note that two polytopes with the same inward normal vectors
define the same toric variety. For example a square and a
rectangle in $\mathbb{Z}^2$ define $\mathbb{P}^1 \times
\mathbb{P}^1$ but they define different Cartier divisors.

A complete fan $\triangle$ and a $T$-invariant Cartier divisor $D=
\sum a_\rho V(\rho)$ defines a polytope,
$$
P_D = \{ u \in \mr \talque  \langle u, v(\rho)\rangle \ge - a_\rho
~ \forall ~ \rho \mathrm{~border~of~} \triangle \}
$$

A toric variety defined from a fan $\triangle$ is normal and it is
projective if and only if $\triangle$ is a fan associated to a
polytope in $\mr$.

The following lemma allows us compute a basis of
$\mathcal{O}(D_P)$.

\begin{lema} \label{le:co}
Let $X_P$ be the toric variety associated to a polytope
$\triangle$. The set $\codh$ of global sections of
$\mathcal{O}(D_P)$ is a finite dimensional $\fq$-vector space with
$\{\chi^u \talque u \in M \cap P \}$ as a basis.
\end{lema}

%%%%%%%%%%%%%%%%%%%%%%%%%%%%%%%%%%%%%%%%%%%%%%%%%%%%%%%%%%%%%%%%%%%%%%
%%%%%%%%%%%%%%%%%%%%%%%%%%%%%%%%%%
%%%%%%%%%%%%%%%%%          %%%%%%%%%%%%%%%           %%%%%%%%%%%%%%%%%
%%%%%%%%%%%%%%%%%%%%%%%%%%%%%%%%%%%%%%%%%%%%%%%%%%%%%%%%%%%%%%%%%%%%%%
%%%%%%%%%%%%%%%%%%%%%%%%%%%%%%%%%%%%%%%%%%%%%%%%%%%%%%%%%%%%%%%%%%%%%%
%%%%%%%%%%%%%%%%%          %%%%%%%%%%%%%%%           %%%%%%%%%%%%%%%%%
%%%%%%%%%%%%%%%%%%%%%%%%%%%%%%%%%%
%%%%%%%%%%%%%%%%%%%%%%%%%%%%%%%%%%%%%%%%%%%%%%%%%%%%%%%%%%%%%%%%%%%%%%

\section{Toric Codes}

Let $P$ be a rational polytope of dimension $r \ge 2$, $X_P$ its
associated refined variety and $D_P$ its associated Cartier
divisor on $X_P$ as in the previous section.

For $t \in T = (\fq^\ast)^r$, the rational functions of $\codh$,
i.e. rational functions $f$ over $X_P$ such that $\mathrm{div}(f)
+ D_P \succeq 0$, can be evaluated at $t$
$$
\begin{array}{ccc}
  \codh & \to & \fq \\
  f & \mapsto & f(t) \\
\end{array}
$$since $f$ is a linear combination of characters $\chi^u$ that can be considered as Laurent monomials (lemma \ref{le:co}). This map is nothing else that the evaluation of a
Laurent polynomial whose monomials have exponents in $P$ in a
point with non-zero coordinates.

We define the toric codes, in the same way as Hansen \cite{ha}.
Evaluating at the $(q-1)^r$ points of $T = (\fq^\ast)^r$ we obtain
the \textbf{toric code $\code_P$ associated to $P$}, which is an
evaluation code in the sense of Goppa \cite{tsvl}. $\code_P$ is
the image of the $\fq$-linear evaluation map given by
$$
\begin{array}{ccc}
  \ev:\codh & \to & \left(\fq \right)^{\# T} \\
  f & \mapsto & (f(t))_{t \in T} \\
\end{array}
$$
Since we evaluate in $\# T$ points, $\code_P$ has \textbf{length}
$n=\# T = (q-1)^r$.

From lemma \ref{le:co}, it follows that $\codh$ is a $\fq$-vector
space of finite dimension with basis $\{\chi^u \talque u \in M
\cap P \}$, therefore a generator system of the code $\code_P$ is
$\{(\chi^u (t))_{t \in T} \talque u \in M \cap P \}$ which is also
a basis of the code if and only if the evaluation map $\ev$ is
injective.

\begin{nota}
D. Joyner in \cite{jo} defines a code for a toric variety coming
from a complete fan, a Cartier divisor and a 1-cycle, Joyner uses
the 1-cycle to evaluate the rational functions in its support.
Then he consider the special case where the 1-cycle has support
$T$ and he call this codes \textbf{standard toric codes}. As we
have seen in the previous section a complete fan and a Cartier
divisor is the same data as a polytope $P$. A polytope $P$
determines the fan $\triangle_P$, the toric variety $X_P$ and the
Cartier Divisor $D_P$. Therefore the toric codes defined here
which are the same as Hansen's construction \cite{ha} are as
general as the standard toric code (\cite[definition 4.5]{jo}) of
Joyner \cite{jo}.
\end{nota}

The following lemma is used to compute the kernel of the
evaluation map and the dimension of the code is given in theorem
\ref{th:nu}.

\begin{lema} \label{le:ce}
Let $P$ be a polytope such that $P \cap M$ is contained in\\$H =
\{ 0, \ldots, q-2 \} \times \cdots \times \{ 0, \ldots, q-2 \}
\subset M$. Let
$$
f=\sum_{u \in P \cap M} \lambda_u \chi^u, \; \; \; \; \lambda_u
\in \fq
$$
Then $(f(t))_{t\in T} = (0)_{t\in T}$ ($f \in \ker (\ev)$ for some
$D$) if and only if $\lambda_u = 0,\\\forall \;  u \in P \cap M$.
\end{lema}

\pf

Let $f= \sum_{u \in P \cap M} \lambda_u \chi^{u}$, we can write
$f$ as
$$
f(t_1, \ldots, t_r) = \sum_{0 \le u_1, \ldots, u_r \le q-2}
\lambda_{u_1, \ldots, u_r} t_1^{u_1} \cdots t_r^{u_r} \in \fq
[t_1, \ldots, t_r]
$$ with $\lambda_{u_1, \ldots, u_r} \in \fq$. We shall see that $f = 0$.

We prove the result by induction in the number of variables. If
$r=1$, $f = \sum_{0 \le u_1 \le q-2} \lambda_{u_1} t_1^{u_1}$,
since $f$ vanish in all $\fq^\ast$ it belongs to the ideal
generated by $t_1^{q-1} -1$, therefore $f = 0$ (by degree
considerations).

Assume that the result holds up to $r-1$ variables. Let $t_1,
\ldots, t_{r-1} \in \fq^\ast$ then $$f(t_1, \ldots, t_{r-1},t_r) =
g_{q-2}(t_1, \ldots, t_{r-1}) t_r^{q-2} + \cdots + g_{1}(t_1,
\ldots, t_{r-1})t_r + g_{o}(t_1, \ldots, t_{r-1}) $$ with
$g_i(t_1, \ldots, t_{r-1}) \in \fq [t_1, ..., t_{r-1}]$.

One has that  $f(t_1, \ldots, t_{r-1},t_r) \in \fq[t_r]$ vanish
for all $t_r \in \fq^\ast$. Therefore $f$ belongs to the ideal
generated by $t_r^{q-1} -1$, then $f=0$ (by degree
considerations). Hence $g_{i} = 0$ for all $i=1, \ldots, q-2$ and
we can apply the induction hypothesis to $g_i$ and we obtain
$f=0$. \qed

The following theorem allows us to compute the kernel of the
evaluation map and a basis of the code (and therefore its
dimension).

\begin{teorema}\label{th:nu}

Let $P$ be a polytope and $\code_P$ be its associated toric code.

For all $u \in P \cap M$ we write $u= c_u + b_u$ where $c_u \in H
= \{ 0, \ldots, q-2 \} \times \cdots \times \{ 0, \ldots, q-2 \}
\subset M$, and $b_u \in ((q-1)\zz)^r$. Let $\overline{P}$ be the
set, $\overline{P} = \{ c_u \talque u \in P \} \subset M$.

\vspace{0.4cm}

One has that,

\begin{itemize}

\item[(1)] The kernel of the evaluation map $\ev$ is the
$\fq$-vector space generated by $$\{ \chi^{u} - \chi^{u'} \talque
u, u' \in P \cap M, \; c_u = c_{u'} \}$$

\item[(2)] A basis of the code $\code_P$ is
$$
\{ (\chi^{c_u}(t))_{t\in T}  \talque u \in P \cap M \} = \{
(\chi^{u}(t))_{t\in T}  \talque u \in \overline{P}  \}
$$

and therefore the \textbf{dimension} of $\code_P$ $$k = \# \{ c_u
\talque u \in P \cap M \} = \# \overline{P}$$
\end{itemize}
\end{teorema} \pf \begin{itemize}

\item[(1)] Let $u, u' \in P \cap M$ such that  $c_u = c_{u'}$.
Then $\ev(\chi^u) = \ev (\chi^{u'})$ and one has that
$\ev(\chi^{u}-\chi^{u'}) \in \ker(\ev)$.

On the other hand let $f \in \codh $, with $\ev(f)=0$.

$$
f =  \sum_{u \in P \cap M} \lambda_u \chi^{u} = \sum_{u \in P \cap
M} \lambda_u (\chi^{u} - \chi^{c_u}) + \sum_{u \in P \cap M}
\lambda_u \chi^{c_u}
$$

One has for all $t \in T$
$$
\underbrace{f(t)}_{=0} =  \underbrace{\sum_{u \in P \cap M}
\lambda_u (\chi^{u}(t) - \chi^{c_u}(t))}_{=0} + \sum_{u \in P \cap
M} \lambda_u \chi^{c_u}(t)
$$

Then $\sum_{u \in P \cap M} \lambda_u  \chi^{c_u}(t) = 0$ for all
$t \in T$, and applying the lemma \ref{le:ce} ($c_u \in H \;
\forall \; u$) one has that $\sum_{u \in P \cap M} \lambda_u
\chi^{c_u}$ is the zero function. Then $f$ belongs to the vector
space generated by $\{ \chi^{u} - \chi^{u'} \talque u, u' \in P
\cap M, \; c_u = c_{u'} \}$.

\vspace{0.5cm}

\item[(2)]

Let $f \in \codh$, and let $t\in T$,
$$
f(t) =  \sum_{u \in P \cap M} \lambda_u \chi^{u}(t) = \sum_{u \in
P \cap M} \lambda_u \chi^{c_u + b_u}(t) = \sum_{u \in P \cap M}
\lambda_u \chi^{c_u} (t)
$$

Therefore $(f(t))_{t \in T} \in \{ (\chi^{c_u}(t))_{t\in T}
\talque u \in P \cap M \}$.

And moreover $\{ (\chi^{c_u}(t))_{t\in T}  \talque u \in P \cap M
\}$ is linear independent set by the lemma \ref{le:ce} ($c_u \in H
\; \forall \; u$).\qed

\end{itemize}

Two polytopes $P$, $P'$ such that $\overline{P} = \overline{P'}$
have the same associated toric code ($\code_P = \code_{P'}$).
Computing $\chi^{c_u}$ is the same as computing the class of
$\chi^u$ in $\fq[X_1, \ldots, X_r]/ J$, where
$J=(X_1^{q-1}-1,\ldots,X_r^{q-1}-1)$. In \cite{diguva} it is
proven that a toric code of dimension 2 is multicyclic,
considering the class of $\chi^u$ in $\fq[X_1, \ldots, X_r]/ J$
one can see that $\code_P$ is multicyclic for arbitrary dimension.

We say that a polytope $P$ verifies the \textbf{injectivity
restriction} if for all $u, u' \in P \cap M, u \neq u'$ one has
that $c_u \neq c_{u'}$. Using the above theorem, $P$ verifies the
injectivity restriction if and only if the evaluation map $\ev$ is
injective and $\code_P$ has therefore dimension $k = \# (P \cap
M)$, that is the number of rational points in the polytope. In
\cite{ha,ha3} Hansen restricts the size of the polytopes in order
to make the evaluation map injective, by considering the minimal
distance bound. The dimension of the code is therefore the number
of rational points of the polytope.

A discussion of recent algorithms to compute the number of lattice
points in a polytope may be found in \cite{de}. For $r=2$ one has
Pick's formula \cite{fu} to compute the number of lattice points:

\begin{lema}

\label{le:pi} Let $P$ be a plane polytope. Then
$$
\# (P \cap M) = \mathrm{vol}_2 (P) +
\frac{\mathrm{Perimeter}(P)}{2} +1
$$where $\mathrm{vol}_2$ is the Lebesgue volume.
\end{lema}

%%%%%%%%%%%%%%%%%%%%%%%%%%%%%%%%%%%%%%%%%%%%%%%%%%%%%%%%%%%%%%%%%%%%%%
%%%%%%%%%%%%%%%%%%%%%%%%%%%%%%%%%%
%%%%%%%%%%%%%%%%%          %%%%%%%%%%%%%%%           %%%%%%%%%%%%%%%%%
%%%%%%%%%%%%%%%%%%%%%%%%%%%%%%%%%%%%%%%%%%%%%%%%%%%%%%%%%%%%%%%%%%%%%%
%%%%%%%%%%%%%%%%%%%%%%%%%%%%%%%%%%%%%%%%%%%%%%%%%%%%%%%%%%%%%%%%%%%%%%
%%%%%%%%%%%%%%%%%          %%%%%%%%%%%%%%%           %%%%%%%%%%%%%%%%%
%%%%%%%%%%%%%%%%%%%%%%%%%%%%%%%%%%
%%%%%%%%%%%%%%%%%%%%%%%%%%%%%%%%%%%%%%%%%%%%%%%%%%%%%%%%%%%%%%%%%%%%%%

\section{Estimates for the minimum distance} \label{se:di}

Finally in order to compute the parameters of this family of codes
we compute the minimum distance. We use the same techniques of
\cite{ha} for dimension 2, and compute the intersection numbers
using mixed volumes. We also extend this computations to arbitrary
dimension. In order to compute the \textbf{minimum distance} $d$
of the linear code $\code_P$ we should compute the minimum weight
of a non-zero word, i.e. the maximum number of zeros of a function
$f$ in $\codh \setminus \{ 0 \}$ in $T$. We solve this problem
using intersection theory.

Let $u_1=(1,0, \ldots,0), u_2=(0,1,0,\ldots,0), \ldots, u_r=(0,
\ldots,0,1)$. Each $\fq$-rational point of $T$ is contained in one
of the $(q-1)^{r-1}$ lines
$$
C_{\eta_1,\ldots,\eta_{r-1}} = \mathrm{Z} ( \{ \chi^{u_i} - \eta_i
: i=1,\ldots, r-1  \} ), ~ ~  ~  ~ ~ \eta_i \in \fq^\ast ~ \forall
i
$$

Let $f \in \codh \setminus \{ 0\}$. Assume that $f$ is identically
zero in $a$ of the lines, and denote by $A$ the set of subindexes
of the $a$ lines where $f$ vanish.

Following \cite[proposition 3.2]{hn}, in the other lines the
number of zeros is given by the intersection number of a Cartier
divisor with a 1-cycle, the integer $D_P \cdot
C_{\eta_1,\ldots,\eta_{r-1}}$. Therefore the number of zeros of
$f$ in $T$ is bounded by

$$
a (q-1) + \sum_{\eta_i \in \fq^\ast,\\ (\eta_1, \ldots,
\eta_{r-1}) \notin A} (D_P \cdot C_{\eta_1,\ldots,\eta_{r-1}} )
$$

In order to compute the maximum number of zeros of $f$ one has to
compute the intersection number of the Cartier divisor and the
1-cycle and bound the number of lines where $f$ is $0$.

Following \cite{fu2}   $D_P \cdot C_{\eta_1,\ldots,\eta_{r-1}} =
D_P \cdot C $ for any $C$ defined above. Therefore the number of
zeros of $f$ is bounded by
$$
a (q-1) + ((q-1)^{r-1} -a) (D_P \cdot C )
$$and the minimum distance is bounded by
$$
d(\code_P) \ge n - ( a(q-1) + ((q-1)^{r-1} -a)(D_P \cdot C))
$$

One has that
$$
D_P \cdot C  = D_P \cdot (\div (\chi^{u_1}) )_0 \cdot \cdots \cdot
(\div (\chi^{u_{r-1}}) )_0
$$and following \cite{fu} one see that this intersection number is the mixed volume of the associated polytopes$$
r!V_r (P, P_{(\div (\chi^{u_1}) )_0}, \cdots, P_{(\div
(\chi^{u_{r-1}}) )_0} )
$$The \textbf{mixed volume} $V_r$ of $r$ polytopes $P_1, \cdots, P_r$ is
$$
V_r (P_1, \ldots, P_r) = \frac{1}{r!}\sum_{j=1}^r (-1)^{r-j}
\sum_{1 \le i_1 < \cdots < i_j \le r} \mathrm{Vol}_r (P_{i_1} +
\cdots + P_{i_j})
$$where $\mathrm{Vol}_r$ is the Lebesgue volume. An algorithm to compute
the Lebesgue volume of a polytope may be found in \cite{ba}.
Moreover under certain hypothesis the mixed volume can be computed
directly \cite{ko}.

Let $f \in \codh$, since $\code_P = \code_{\overline{P}}$ we
assume without loss of generality that $\deg_{t_i} f \le q-2$.
$$
f(t_1, \ldots, t_r ) = f_0 (t_1, \ldots, t_{r-1}) + f_1 (t_1,
\ldots, t_{r-1}) t_r + \cdots + f_{q-2} (t_1, \ldots,
t_{r-1})t_r^{q-2}
$$let $C_{\eta_1,\ldots,\eta_{r-1}}$ be a line where $f$ vanish, $f(\eta_1,
\ldots , \eta_{r-1}, t_r) \in \fq[t_r]$, and\\$\deg f(\eta_1,
\ldots , \eta_{r-1}, t_r) < t_r^{q-1}$ therefore since $f(\eta_1,
\ldots , \eta_{r-1}, t_r) = 0 ~ \forall ~ t_r \in \fq^\ast$ one
has that $f_i (\eta_1, \ldots, \eta_{r-1} ) = 0 ~ \forall ~ i$.

The number $a$ is less than or equal to the maximum number of
zeros of a non zero function $f \in \mathrm{H}^0
(X_P',\mathcal{O}(D_P'))$ where $P'$ is the $r$-projection of the
polytope $P$. This can be repeated until we reach dimension 2.
\vspace{0.5cm}

For a \textbf{plane polytope} we compute the minimum distance as
in \cite{ha3}.

Let us consider $P$ a plane polytope and let us bound the minimum
distance. In dimension 2 we can improve the previous computation.
Let $f \in \codh \setminus \{ 0 \}$, and let us assume that $f$ is
identically 0 in $a$ lines. Therefore following \cite[proposition
3.2]{hn} in the other $(q-1-a)$ lines the maximum number of zeros
is $D_P \cdot \div (\chi^{u_1})$.

In dimension 2 a 1-cycle is a Weil divisor and since $f$ vanish in
$a$ of the previous lines one has that
$$
\div (f) + D_P - a ( \div (\chi^{u_1}) )_0 \succeq 0
$$

Or equivalently, $f \in \mathrm{H}^0 (X_P,\mathcal{O}(D_P - a (
\div (\chi^{u_1}) )_0 ))$, and the maximum number of zeros of $f$
in the other $(q-1-a)$ lines is $D_P - a ( \div (\chi^{u_1}) )_0
\cdot ( \div (\chi^{u_1}) )_0$, which is smaller than or equal to
the previous one. This will probably allow us to to give a sharper
bound.

From lemma \ref{le:co} one has that
$$
a \le \max \{ u_2 - u'_2 \talque u_1 = u'_1, (u_1,u_2) \in P,
(u'_1,u'_2) \in P \}
$$

Finally we compute the intersection number of the two Cartier
divisors just in the same way as for $r >2$, using the mixed
volume of the associated polytopes:
$$
D_P - a ( \div (\chi^{u_1}) )_0 \cdot ( \div (\chi^{u_1}) )_0  = 2
V_2 (P_{D_P - a ( \div (\chi^{u_1}) )_0}    , P_{(\div
(\chi^{u_1}) )_0})
$$

\begin{nota}
For a polytope $P$ large enough one can obtain a trivial bound of
the minimum distance, which is not the case when the injectivity
restriction is satisfied. For instance if we consider a rectangle
$P$ with a basis of length greater than $q-1$ we obtain a negative
bound of the minimum distance. Another possibility may be to apply
the above computations to $\overline{P}$ to obtain a non trivial
bound but unfortunately $\overline{P}$ is not in general a convex
polytope. This is similar to the situation for an AG-code $L(D,G)$
when $n \le 2g -2 \deg (G)$ \cite{tsvl}.
\end{nota}

The following proposition gives an \textbf{upper bound} of the
minimum distance and in particular it may be used to check if the
previous bound is sharp. This result extends the computations of
\cite{jo,ha3}.

\begin{proposicion} \label{pr:up}

Let $P$ be a polytope and $\code_P$ its associated linear code.

Let $u \in M$ and $Q$ be $\{0, 1, \ldots, l_1 \} \times \cdots
\times \{ 0, 1, \ldots, l_r \} \subset M$,  where $0 \le l_i \le
q-2$ (some $l_i$ can be equal to zero), if $\overline{u + Q}$ is
contained into the set $\overline{P}$, (where $u=c_u + b_u$, $c_u
\in H$, $b_u \in ((q-1)\mathbb{Z})^r$, $\overline{P}=\{ c_u
\talque u \in P \cap M \}$ as in theorem \ref{th:nu}) then
$$
d \le n - \sum_{j=1}^{r} (-1)^{j+1} \sum_{i_1 < \cdots < i_j}
l_{i_1} \cdots l_{i_j} (q-1)^{r-j}
$$
\end{proposicion}

\pf

 Let $a^i_1, a^i_2, \ldots, a^i_{l_i} \in  \fq^\ast$ be
pairwise different elements for $i=1, \ldots, r$.

Let $f (t_1, \ldots, t_r) = t_1^{u_1} \cdots t_r^{u_r} \prod (t_i
- a^i_1) \cdots (t_i - a^i_{l_i})$. The number of zeros of $f$ in
$T$ is equal to $\sum_{j=1}^{r} (-1)^{j+1} \sum_{i_1 < \cdots <
i_j} l_{i_1} \cdots l_{i_j} (q-1)^{r-j}$ (by the
inclusion-exclusion principle).

Since $f$ is a linear combination of monomials with exponents in
$(u + Q) \cap M$ and $\overline{u + Q} \subset \overline{P}$ one
has that for each monomial $\chi^{c_u}$ in $f$ there exists $b_u
\in ((q-1)\mathbb{Z})^r$ such that $\chi^{c_u + b_u} \in \codh$,
and both polynomials take the same values in $T$. Proceeding in
the same way with all the monomials of $f$ one obtains a function
$f'$ such that $f'(t)=f(t)$, $\forall t \in T$ and $f' \in \codh$.
Therefore an upper bound for the minimum distance is

$$
d \le n - \sum_{j=1}^{r} (-1)^{j+1} \sum_{i_1 < \cdots < i_j}
l_{i_1} \cdots l_{i_j} (q-1)^{r-j} \qed
$$

%%%%%%%%%%%%%%%%%%%%%%%%%%%%%%%%%%%%%%%%%%%%%%%%%%%%%%%%%%%%%%%%%%%%%%
%%%%%%%%%%%%%%%%%%%%%%%%%%%%%%%%%%
%%%%%%%%%%%%%%%%%          %%%%%%%%%%%%%%%           %%%%%%%%%%%%%%%%%
%%%%%%%%%%%%%%%%%%%%%%%%%%%%%%%%%%%%%%%%%%%%%%%%%%%%%%%%%%%%%%%%%%%%%%
%%%%%%%%%%%%%%%%%%%%%%%%%%%%%%%%%%%%%%%%%%%%%%%%%%%%%%%%%%%%%%%%%%%%%%
%%%%%%%%%%%%%%%%%          %%%%%%%%%%%%%%%           %%%%%%%%%%%%%%%%%
%%%%%%%%%%%%%%%%%%%%%%%%%%%%%%%%%%
%%%%%%%%%%%%%%%%%%%%%%%%%%%%%%%%%%%%%%%%%%%%%%%%%%%%%%%%%%%%%%%%%%%%%%

\section{Examples}

We consider two examples. We first ilustrate the computations of
the parameters for a sequence of polytopes $(P_r)_{r \ge 2}$ with
$\dim (P_r) = r$ and when the $r$-projection of $P_r$ is
$P_{r-1}$. The second example shows that the bound of the minimum
distance, using intersection theory, does not equal to the upper
bound of the proposition \ref{pr:up}

\begin{ejemplo}

Let $P_2$ be the plane polytope of vertices $(0,0)$, $(b_1,0)$,
$(b_1,b_2)$, $(0,b_2)$ with $b_1, b_2 < q - 1$. This is the code
of \cite[proposition 3.2]{ha}.

The fan $\triangle_{P_2}$ associated to $P_2$ is generated by
cones where the edges are generated by $v(\rho_1)) =(1,0)$,
$v(\rho_2) =(0,1)$, $v(\rho_3) =(-1,0)$ and $v(\rho_4) =(0,-1)$.
The toric variety $X_{P_2}$ is non-singular.

$$
P_2 = \bigcap_{i=1}^{4} \{ \langle u, \rho_i \rangle \ge - a_i \}
$$where $a_1= 0 $, $a_2 = 0 $, $a_3 = b_1 $, $a_4 = b_2 $. Therefore $D_P = \sum a_i V(\rho_i) = b_1 V(\rho_3) + b_2 V(\rho_4)$.

Since $P_2$ is a plane polytope the code $\code_{P_2}$ has length
$n=(q-1)^2$. The evaluation map $\ev$ is injective since $b_1, b_2
< q-1$ and $P_2$ verifies the injectivity restriction \ref{th:nu}.
Therefore one has that the dimension of $\code_{P_2}$ is

$$
k= \dim \mathrm{H}^0 (X_{P_2},\mathcal{O}(D_{P_2})) = \# P_2 \cap
M = (b_1 + 1)(b_2 + 1)
$$

From section \ref{se:di} we get that the maximum number of zeros
of a function $f$ in $\mathrm{H}^0 (X_{P_2},\mathcal{O}(D_{P_2}))$
is smaller than or equal to
$$
a(q-1) + (q-1-a) (D_{P_2} - a(\div ( \chi^{u_1} ))_0  \cdot (\div
( \chi^{u_1} ))_0)
$$where $a \le b_1$.

One has that $\div(\chi^{u_1}) = \sum \langle u_1, v(\rho_i)
\rangle V(\rho_i) = V(\rho_1) - V(\rho_3)$. Therefore\\$(\div (
\chi^{u_1} ))_0 = V(\rho_1)$.

$$
D_{P_2}-a(\div ( \chi^{u_1} ))_0  \cdot (\div ( \chi^{u_1} ))_0 =
2 V_2 ( P_{D_{P_2}-a(\div ( \chi^{u_1} ))_0}, P_{(\div (
\chi^{u_1} ))_0} )=
$$ $$
\vold ( P_{D_{P_2}-a(\div ( \chi^{u_1} ))_0} + P_{(\div (
\chi^{u_1} ))_0} ) - \vold ( P_{D_{P_2}-a(\div ( \chi^{u_1} ))_0})
- \vold ( P_{(\div ( \chi^{u_1} ))_0} )=
$$
$$
((b_1-a+1)b_2) - ((b_1 -a)b_2) -(0) = b_2
$$

Because
\begin{itemize}

\item $P_{D_{P_2}-a(\div ( \chi^{u_1} ))_0} + P_{(\div (
\chi^{u_1} ))_0}$ is the polytope of vertices $(a-1,0)$,
$(b_1,0)$, $(b_1,b_2)$ and $(a-1,b_2)$. \item$P_{D_{P_2}-a(\div (
\chi^{u_1} ))_0}$ is the polytope of vertices $(a,0)$, $(b_1,0)$,
$(b_1,b_2)$ and $(a,b_2)$. \item $P_{(\div ( \chi^{u_1} ))_0}$ is
the polytope of vertices $(-1,0)$ and $(0,0)$.

\end{itemize}

Therefore the maximum number of zeros of $f\in \mathrm{H}^0
(X_{P_2},\mathcal{O}(D_{P_2}))$ is bounded by
$$
a (q-1-b_2) + (q-1)b_2 \le b_1(q-1-b_2) + (q-1)b_2
$$

and the minimum distance is bounded by
$$
d \ge (q-1)^2 - (b_1(q-1-b_2) + (q-1)b_2) = (q-1-b_1)(q-1-b_2)
$$

We then apply proposition \ref{pr:up}, with $u=0$ and $l_1 = b_1$,
$l_2 =b_2$. $u + Q \subset P_2$, then indeed $u +Q = P_2$, and we
obtain
$$d \le (q-1)^2 - b_1(q-1) -b_2 (q-1) + b_1b_2 = (q-1-b_1)(q-1-b_2)$$
And therefore $d = (q-1-b_1)(q-1-b_2)$

%%%%

%%%5

Let $P_3$ be the 3 dimensional polytope of vertices $(0,0,0)$,
$(b_1,0,0)$, $(b_1,b_2,0)$, $(0,b_2,0)$, $(0,0,b_3)$,
$(b_1,0,b_3)$, $(b_1,b_2,b_3)$, $(0,b_2,b_3)$ with $b_1, b_2, b_3
< q - 1$.

The fan $\triangle_{P_3}$ associated to $P_3$ is generated by
cones with edges generated by $v(\rho_1)) =(1,0,0)$, $v(\rho_2)
=(-1,0,0)$, $v(\rho_3) =(0,1,0)$, $v(\rho_4) =(0,-1,0)$,
$v(\rho_5) =(0,0,1)$, $v(\rho_6) =(0,0,-1)$. The toric variety
$X_{P_3}$ is non-singular.

$$
P_3 = \bigcap_{i=1}^{6} \{ \langle u, \rho_i \rangle \ge - a_i \}
$$where $a_1= 0 $, $a_2 = b_1 $, $a_3 = 0 $, $a_4 = b_2 $, $a_5 = 0 $, $a_6 = b_3 $. Therefore $D_P = \sum a_i V(\rho_i) = b_1 V(\rho_2) + b_2 V(\rho_4) +b_3 V(\rho_6)$.

Since $P_3$ is a 3 dimensional polytope the code $\code_{P_3}$ has
length $n=(q-1)^3$. The evaluation map $\ev$ is injective since
$b_1, b_2, b_3 < q-1$ and $P_3$ verifies the injectivity
restriction \ref{th:nu}. Therefore one has that the dimension of
$\code_{P_3}$ is
$$
k= \dim \mathrm{H}^0 (X_{P_3},\mathcal{O}(D_{P_3})) = \# P_3 \cap
M = (b_1 + 1)(b_2 + 1)(b_3 +1)
$$

From section \ref{se:di} the maximum number of zeros of a function
$f\in\mathrm{H}^0 (X_{P_3},\mathcal{O}(D_{P_3}))$ is smaller than
or equal to
$$
a(q-1) + ((q-1)^2 -a) (D_{P_3} \cdot C)
$$where $C= \mathrm{Z} ( \{ \chi^{u_1} , \chi^{u_2}  \})$ and $a$ is
smaller than or equal to the maximum number of zeros of a function
defined by the 3-projection of $P_3$, i.e. $P_2$. Therefore $a \le
b_1(q-1-b_2) + (q-1)b_2$.

One has that $\div(\chi^{u_1}) = \sum \langle u_1, v(\rho_i)
\rangle V(\rho_i) = V(\rho_1) - V(\rho_2)$.
Therefore\\$(\div(\chi^{u_1} ))_0 = V(\rho_1)$. $\div(\chi^{u_2})
= \sum \langle u_1, v(\rho_i) \rangle V(\rho_i) = V(\rho_3) -
V(\rho_4)$. Therefore $(\div ( \chi^{u_2} ))_0 = V(\rho_3)$.

$$
D_{P_3} \cdot C = D_{P_3} \cdot (\div(\chi^{u_1}))_0 \cdot
(\div(\chi^{u_2}))_0 = 3! V_3
(P,P_{(\div(\chi^{u_1}))_0},P_{(\div(\chi^{u_2}))_0})=
$$
$$
\volt ( P_3 + P_{(\div ( \chi^{u_1} ))_0} + P_{(\div ( \chi^{u_2}
))_0} ) - \volt ( P_3 + P_{(\div ( \chi^{u_1} ))_0}) - \volt ( P_3
+ P_{(\div ( \chi^{u_2} ))_0}) -$$$$- \volt ( P_{(\div (
\chi^{u_1} ))_0} + P_{(\div ( \chi^{u_2} ))_0}) + \volt(P_3) +
\volt(P_{(\div ( \chi^{u_1} ))_0}) + \volt(P_{(\div ( \chi^{u_2}
))_0})=
$$
$$
((b_1+1)(b_2+1)(b_3)) -((b_1+1)b_2 b_3 ) -( b_1(b_2 +1)b_3) - (0)
+ (b_1 b_2 b_3) + (0) +(0) = b_3
$$

Because
\begin{itemize}

\item $P_3 + P_{(\div ( \chi^{u_1} ))_0} + P_{(\div ( \chi^{u_2}
))_0}$   is the polytope of vertices  $(-1,-1,0)$,\\$(b_1,-1,0)$,
$(b_1,b_2,0)$, $(-1,b_2,0)$, $(-1,-1,b_3)$, $(b_1,-1,b_3)$,
$(b_1,b_2,b_3)$ and $(-1,b_2,b_3)$. \item $P_3 + P_{(\div (
\chi^{u_1} ))_0}$   is the polytope of vertices  $(-1,0,0)$,
$(b_1,0,0)$, $(b_1,b_2,0)$, $(-1,b_2,0)$, $(-1,0,b_3)$,
$(b_1,0,b_3)$, $(b_1,b_2,b_3)$ and $(-1,b_2,b_3)$. \item  $P_3 +
P_{(\div ( \chi^{u_2} ))_0}$  is the polytope of vertices
$(0,-1,0)$, $(b_1,-1,0)$, $(b_1,b_2,0)$, $(0,b_2,0)$,
$(0,-1,b_3)$, $(b_1,-1,b_3)$, $(b_1,b_2,b_3)$ and $(0,b_2,b_3)$.
\item $P_{(\div ( \chi^{u_1} ))_0} + P_{(\div ( \chi^{u_2} ))_0}$
is the polytope of vertices $(0,0,0)$, $(-1,0,0)$, $(-1,-1,0)$ and
$(0,-1,0)$. \item  $P_3$  is the polytope of vertices  $(0,0,0)$,
$(b_1,0,0)$, $(b_1,b_2,0)$, $(0,b_2,0)$, $(0,0,b_3)$,
$(b_1,0,b_3)$, $(b_1,b_2,b_3)$ and $(0,b_2,b_3)$. \item  $P_{(\div
( \chi^{u_1} ))_0}$  is the polytope of vertices $(-1,0,0)$ and
$(0,0,0)$. \item  $P_{(\div ( \chi^{u_2} ))_0}$  is the polytope
of vertices  $(0,-1,0)$ and $(0,0,0)$

\end{itemize}

Therefore the maximum number of zeros of $f\in \mathrm{H}^0
(X_{P_3},\mathcal{O}(D_{P_3}))$ is bounded by
$$
a (q-1- b_3) + (q-1)^2 b_3 \le  (b_1(q-1-b_2) + (q-1)b_2)(q-1-b_3)
+ (q-1)^2 b_3
$$and the minimum distance is bounded by
$$
d \ge n - ( (b_1(q-1-b_2) + (q-1)b_2) (q-1-b_3) + (q-1)^2 b_3) =
(q-1-b_1)(q-1-b_2)(q-1-b_3)
$$

We then apply proposition \ref{pr:up}, with $u=0$ and $l_1 = b_1$,
$l_2 =b_2$, $l_3=b_3$. $u + Q \subset P_3$, then indeed $u +Q =
P_3$, and we obtain
$$d \le (q-1-b_1)(q-1-b_2) (q-1-b_3)$$

And therefore $d = (q-1-b_1)(q-1-b_2)(q-1-b_3)$.

Computing the lower and upper bound of the minimum distance for an
hypercube $P_r$ of dimension $r$ with sides $b_1, \ldots, b_r <
q-1$ one obtain for all $r\ge2$ that its minimum distance $d_r$ is
equal to
$$
d_2 = (q-1-b_1)(q-1-b_2)
$$
$$
d_r = (q-1)^r - ( (q-1)^{r-1} - d_{r-1})(q-1-b_r) -
b_r(q-1)^{r-1}, ~ ~ ~ \forall ~ r \ge 3
$$
one can easily see (by induction on $r$) that it is equal to
$$
d_r = (q-1-b_1) \cdots (q-1-b_r)
$$

Therefore, the code $\code_{P_r}$ associated to the hypercube of
sides $b_1, \ldots, b_r$ has parameters $[(q-1)^r,\prod(b_i + 1),
\prod(q-1-b_i)]$. \cite{li} also consider this example. There the
distance is computed using Vandermonde determinants.

\end{ejemplo}

%%%%%%%%

\vspace{0.5cm}

%%%%%%%

In Hansen's examples \cite{ha3} for plane polytopes and also in
the previous example the lower bound of the minimum distance,
using intersection theory, equals to the upper bound of the
proposition \ref{pr:up}. One could think that the previous bound
is always sharp, the following example shows that this bound is
not always sharp.

\begin{ejemplo}
Let $P$ be the plane polytope of vertices $(0,0)$, $(b,0)$,
$(2b,b)$, $(2b,2b)$, $(b,2b)$, $(0,b)$ with $b < q - 1$.

The fan $\triangle_P$ associated to $P$ is generated by cones with
edges generated by $v(\rho_1)) =(1,0)$, $v(\rho_2) =(0,1)$,
$v(\rho_3) =(-1,1)$, $v(\rho_4) =(-1,0)$, $v(\rho_5) =(0,-1)$,
$v(\rho_6) =(1,-1)$. The toric variety $X_P$ is non-singular.

$$
P = \bigcap_{i=1}^{6} \{ \langle u, v(\rho_i) \rangle \ge - a_i \}
$$where $a_1= 0 $, $a_2 = 0$, $a_3 = b $, $a_4 = 2b  $, $a_5 = 2b $, $a_6 =
b$. Therefore $D_P = \sum a_i V(\rho_i) = b V(\rho_3) + 2b
V(\rho_4) + 2b V(\rho_5) +V(\rho_6)$.

Since $P$ is a plane polytope the code $\code_P$ has length
$n=(q-1)^2$. The evaluation map $\ev$ is injective since $b < q-1$
and $P$ verifies the injectivity restriction. Therefore one has
that dimension of $\code_P$ is

$$
k= \dim \codh = \mathrm{vol}_2 (P) +
\frac{\mathrm{Perimeter}(P)}{2} +1 = 3b^2 + 3b +1
$$

From section \ref{se:di} the maximum number of zeros of a function
$f \in \codh$ is smaller than or equal to
$$
a(q-1) + (q-1-a) (D_P-a(\div ( \chi^{u_1} ))_0  \cdot (\div (
\chi^{u_1} ))_0)
$$where $a \le 2b$.

One has that $\div(\chi^{u_1}) = \sum \langle u_1, v(\rho_i)
\rangle V(\rho_i) = V(\rho_1) - V(\rho_3) - V(\rho_4) +
V(\rho_6$). Therefore $(\div ( \chi^{u_1} ))_0 = V(\rho_1) +
V(\rho_6)$.

$$
D_P-a(\div ( \chi^{u_1} ))_0  \cdot (\div ( \chi^{u_1} ))_0 = 2
V_2 ( P_{D_P-a(\div ( \chi^{u_1} ))_0}, P_{(\div ( \chi^{u_1}
))_0} )=
$$ $$
\vold ( P_{D_P-a(\div ( \chi^{u_1} ))_0} + P_{(\div ( \chi^{u_1}
))_0} ) - \vold ( P_{D_P-a(\div ( \chi^{u_1} ))_0}) - \vold (
P_{(\div ( \chi^{u_1} ))_0} )=
$$ $$
(3b^2 -2ab +2b) - (3b^2 -2ab) -(0) = 2b
$$

Because
\begin{itemize}

\item $P_{D_P-a(\div ( \chi^{u_1} ))_0} + P_{(\div ( \chi^{u_1}
))_0}$ is the polytope of vertices $(a-1,0)$, $(b,0)$, $(2b,b)$,
$(2b,2b)$, $(b+a-1,2b)$ and $(a-1,b-a)$. \item$P_{D_P-a(\div (
\chi^{u_1} ))_0}$ is the polytope of vertices $(a,0)$, $(b,0)$,
$(2b,b)$, $(2b,2b)$,\\$(b+a,2b)$ and $(a,b-a)$. \item $P_{(\div (
\chi^{u_1} ))_0}$ is the polytope of vertices $(-1,0)$ and
$(0,0)$.

\end{itemize}

Therefore the maximum number of zeros of $f\in \codh$ is bounded
by
$$
a (q-1-2b) + (q-1)2b \le 2b(q-1-2b) + (q-1)2b = 4b(q-1) - 4b^2
$$and the minimum distance is bounded by

$$
d \ge n - (4b(q-1) - 4b^2) = (q-1)^2 - 4b(q-1) + 4b^2
$$

As we claimed before in this example the lower bound is different
from the upper bounds. One can apply proposition \ref{pr:up} by
considering a segment of length at most $2b$ and a square of side
at most $b$ inside $P$.

Let $u = (0,b)$ and $Q= \{ 0, 1, \ldots, 2b \} \times \{ 0 \}$,
$u+Q \subset P$. Therefore\\$d \le (q-1)^2 - 2b(q-1)$.

Let $u = (0,0)$ and $Q= \{ 0, 1, \ldots, b \} \times \{ 0, 1,
\ldots, b \}$, $u+Q \subset P$. Therefore $d \le (q-1)^2 -
(2b(q-1) - b^2)$.

Then $(q-1)^2 - 4b(q-1) + 4b^2 < (q-1)^2 - 2b(q-1) < (q-1)^2 -
(2b(q-1) -b^2)$.

\end{ejemplo}

%%%%%%%%%%%%%%%%%%%%%%%%%%%%%%%%%%%%%%%%%%%%%%%%%%%%%%%%%%%%%%%%%%%%%%
%%%%%%%%%%%%%%%%%%%%%%%%%%%%%%%%%%
%%%%%%%%%%%%%%%%%          %%%%%%%%%%%%%%%           %%%%%%%%%%%%%%%%%
%%%%%%%%%%%%%%%%%%%%%%%%%%%%%%%%%%%%%%%%%%%%%%%%%%%%%%%%%%%%%%%%%%%%%%
%%%%%%%%%%%%%%%%%%%%%%%%%%%%%%%%%%%%%%%%%%%%%%%%%%%%%%%%%%%%%%%%%%%%%%
%%%%%%%%%%%%%%%%%          %%%%%%%%%%%%%%%           %%%%%%%%%%%%%%%%%
%%%%%%%%%%%%%%%%%%%%%%%%%%%%%%%%%%
%%%%%%%%%%%%%%%%%%%%%%%%%%%%%%%%%%%%%%%%%%%%%%%%%%%%%%%%%%%%%%%%%%%%%%

\section{Joyner's questions and conjectures}

The \emph{question 3.4} of \cite{jo} asks \emph{under what
conditions (if any) is the map $\ev$ an injection}. Our theorem
\ref{th:nu} answers this question completely for standard toric
codes.

We shall prove that the conjectures 4.2 and 4.3 of \cite{jo} are
not true. As a counterexample we consider a code of the theorem
1.2 of \cite{ha3} and a code of theorem 1.3 of \cite{ha3}
respectively.

\begin{conjetura}
\cite[Conjecture 4.2]{jo}: Let $\code(E,D,X)$ \cite[definition
(5), section 3.1]{jo} be the toric code associated to de 1-cycle
$E$, the $T$-invariant Cartier divisor $D$ and the toric variety
$X$. Let

\begin{itemize}
\item $X$ be a non-singular toric variety of dimension $r$. \item
$n$ be so large that there is an integer $N>1$ such that $2N \volr
(P_D) \le n \le 2N^2 \volr (P_D)$
\end{itemize}

If $q$ is  ``sufficiently large" then any $f \in \mathrm{H}^0
(X,\mathcal{O}(D))$ has no more than $n$ zeros in the rational
points of $X$. Consequently,
$$
d \ge n - 2N \volr (P_D)
$$

Here ``sufficiently large" may depend on $X$, $C$ and $D$ but not
on $f$.

\end{conjetura}

\begin{contraejemplo}
We give a counterexample to the previous conjecture. Let $\code_P$
be the code associated to the plane polytope $P$ of vertices
$(0,0)$, $(1,1)$, $(0,2)$. Following \cite{ha3} $\code_P$ has
length $n = (q-1)^2$ and minimum distance equals $d = (q-1)^2 -
2(q-1)$. The non-singular toric variety $X$ is $X_\triangle$,
where $\triangle$ is the fan generated by cones where the edges
generated by $v(\rho_1)) =(1,0)$, $v(\rho_2) =(-1,1)$, $v(\rho_3)
=(-1,0)$, $v(\rho_4) =(-1,-1)$. $E$ is the formal sum of all the
points of $T$ because $\code_P$ is a standard toric code. We
consider $D=D_P$, that is the Cartier divisor associated to $P$,
$D = V(\rho_3) + V(\rho_4)$ and that $\volr (P_D) = \volr (P) =
1$.

From theorem \ref{th:nu} we know that $q$ ``sufficiently large"
means $q \ge 3$. We claim that the conjecture does not hold for $q
\ge 5$, let $q$ be greater or equal than 5 and $N = q-2$.

$$2N \volr (P_D) \le n \le 2N^2 \volr (P_D) \Leftrightarrow 2(q-2) \le (q-1)^2 \le
2(q-2)^2$$that holds for $q \ge 5$.

The conjecture claims that the minimum distance satisfies

$$
d \ge n - 2N \volr (P_D) = (q-1)^2 - 2 (q-2) > (q-1)^2 - 2(q-1) =
d
$$therefore for $q \ge 5$ the conjecture gives a lower bound strictly
greater than the minimum distance, so this is not true.

\end{contraejemplo}

\begin{conjetura}
\cite[Conjecture 4.3]{jo}: Let $\code(E,D,X)$ \cite[definition
(5), section 3.1]{jo} be the toric code associated to de 1-cycle
$E$, the $T$-invariant Cartier divisor $D$ and the toric variety
$X$. Let

\begin{itemize}
\item $X$ be a non-singular toric variety of dimension $r$. \item
$\psi_D (v) = \min_{u \in P_D \cap M} \langle u,v \rangle$ be
strictly convex \item $\deg (C) > \deg (D^r)$
\end{itemize}

If $q$ is ``sufficiently large" then any $f \in \mathrm{H}^0
(X,\mathcal{O}(D))$ has no more than $n$ zeros in the rational
points of $X$. Consequently,
$$
k \ge \dim \mathrm{H}^0 (X,\mathcal{O}(D)) = \# P_D \cap M
$$
$$
d \ge n - r! (\# P_D \cap M)
$$

Moreover if $n > r! (\# P_D \cap M)$ then $\dim \mathrm{H}^0
(X,\mathcal{O}(D)) = \# P_D \cap M$

\end{conjetura}

\begin{contraejemplo}

We give a counterexample to the previous conjecture. Let $\code_P$
be the code associated to the plane polytope $P$ of vertices
$(0,0)$, $(1,0)$, $(0,1)$. Following \cite{ha3} $\code_P$ has
length $n = (q-1)^2$ and minimum distance equals to $d = (q-1)^2 -
(q-1)$. The non-singular toric variety $X$ is $X_\triangle$, where
$\triangle$ is the fan generated by cones where the edges are
generated by $v(\rho_1)) =(1,0)$, $v(\rho_2) =(0,1)$, $v(\rho_3)
=(-1,-1)$, i.e. $X=\mathbb{P}^2$. $E$ is the formal sum of all the
points of $T$ because $\code$ is a standard toric code. We
consider $D=D_P$, that is the Cartier divisor associated to $P$,
$D = V(\rho_3)$, therefore one has that $\psi_D$ is strictly
convex (see \cite[pag 70]{fu}). One has for $P=P_D$ that $\#P\cap
M = 3$. And $(q-1)^2 = \deg (E)
> \deg (D)=1$

From theorem \ref{th:nu} we know that ``sufficiently large" means
$q \ge 3$. We claim than the conjecture does not hold for $q \ge
8$, let $q$ be greater or equal than 8.

The conjecture claims that the minimum distance satisfies

$$
d \ge n - r! (\# P_D \cap M) = (q-1)^2 - 2\cdot 3 > (q-1)^2 -
(q-1) = d
$$therefore for $q \ge 8$ the conjecture gives a lower bound strictly
greater than the minimum distance, so this is not true.
\end{contraejemplo}

\vspace{1cm}

\textbf{Acknowledgments:} The author wishes to thank T. H\o holdt
and F.J. Monserrat for helpful comments on this paper.

\bibliographystyle{hplain}
\bibliography{toricffaAX}

\end{document}